\documentclass[11pt, twoside, a4paper, english, reqno]{amsart}
\usepackage{amssymb,amsfonts,latexsym,color}
\usepackage{graphics,verbatim}
\usepackage{graphicx}
\newtheorem{theorem}{Theorem}[section]

\newtheorem{lemma}{Lemma}[section]

\newtheorem{remark}{Remark}[section]
\newtheorem{definition}{Definition}[section]

\newcommand{\<}{\left\langle}
\renewcommand{\>}{\right\rangle}

\newcommand{\eps}{\varepsilon}

\newcommand{\be} {\begin{equation}}
\newcommand{\ee} {\end{equation}}
\newcommand{\bea} {\begin{eqnarray}}
\newcommand{\eea} {\end{eqnarray}}
\newcommand{\Bea} {\begin{eqnarray*}}
	\newcommand{\Eea} {\end{eqnarray*}}
\newcommand{\pa} {\partial}

\newcommand{\al} {\alpha}
\newcommand{\ba} {\beta}

\newcommand{\ga} {\gamma}
\newcommand{\Ga} {\Gamma}
\newcommand{\Om} {\Omega}

\newcommand{\De} {\Delta}
\newcommand{\la} {\lambda}

\newcommand{\La} {\Lambda}

\newcommand{\noi} {\noindent}

\newcommand{\R}{\mathbb R}
\newcommand{\N}{\mathbb N}
\newcommand{\Rn}{\mathbb R^N}
\newcommand{\Iom}{\int_{\Omega}}
\newcommand{\deb}{\rightharpoonup}
\catcode`\@=11

\makeatletter \@addtoreset{equation}{section} \makeatother

\begin{document}
	
	\title[On the existence of three non-negative solutions for a $(p,q)$-Laplacian system]{On the existence of three non-negative solutions for $(p,q)$-Laplacian system}

\author{Debangana Mukherjee}
\address{Department of Mathematics and Statistics, Masaryk University, 61137 Brno,  Czech Republic}
\email{mukherjeed@math.muni.cz, \,\, debangana18@gmail.com}
	\author{Tuhina Mukherjee}
\address{Department of Mathematics, National Institute of Technology Warangal, Telangana-506004 }
\email{tulimukh@gmail.com}
	
	\subjclass[2010]{Primary 35R11, 35J20, 49J35,  secondary 47G20, 45G05}
	\keywords{(p,q)-fractional Laplacian, Semilinear elliptic system, Weak solution, Non-negative solutions, Variational methods. }
	\date{}

 \begin{abstract}
The present paper studies the existence of weak solutions for
\begin{equation*}
(\mathcal{P})
\left\{\begin{aligned}
(-\Delta)^{s_1}_{p_1} u  &=\la f_1\,(x,u,v) +g_1(x,u)  \,\mbox{ in }\, \Om, \\
(-\Delta)^{s_2}_{p_2} v  &=\la f_2\,(x,u,v) +g_2(x,v) \,\mbox{ in }\, \Om, \\
u=v &= 0  \,\mbox{in }\, \Rn \setminus \Om, \\
\end{aligned}
\right.
\end{equation*}
where $\Om \subset \Rn$ is a smooth bounded domain with smooth boundary, $s_1,s_2 \in (0,1)$, $1<p_i<\frac{N}{s_i}$, $i=1,2$, $f_i$ and $g_i$ has certain growth assumptions for $i=1,2$. We prove existence of at least three  non negative solutions of $(\mathcal P)$ under restrictive range of $\lambda$ using variational methods. As a consequence, we also conclude that a similar result can be obtained when we consider a more general non local operator $\mathcal L_{\phi_i}$ instead of $(-\Delta)^{s_i}_{p_i}$ in $(\mathcal P)$.
\end{abstract}

\maketitle

\tableofcontents
\section{Introduction}
In the present article, we study the following non-local system of semilinear elliptic equations:
\begin{equation*}
(\mathcal{P})
\left\{\begin{aligned}
(-\Delta)^{s_1}_{p_1} u  &=\la f_1\,(x,u,v) +g_1(x,u)  \,\mbox{ in }\, \Om, \\
(-\Delta)^{s_2}_{p_2} v  &=\la f_2\,(x,u,v) +g_2(x,v) \,\mbox{ in }\, \Om, \\
u=v &= 0  \,\mbox{in }\, \Rn \setminus \Om, \\
\end{aligned}
\right.
\end{equation*}
where $\Om \subset \Rn$ is a smooth bounded domain in with smooth boundary, $s_1,s_2 \in (0,1)$, $1<p_i<\frac{N}{s_i}$, $i=1,2$, the operators $(-\De)^{s_i}_{p_i}$ for $ i=1,2$ are defined as:
\begin{align} \label{frac s_p}
	(-\Delta)^{s_i}_{p_i} u(x)=\lim_{\eps\to 0}\int_{\mathbb{R}^N\setminus B_\eps(x)}\frac{|u(y)-u(x)|^{p_i-2}(u(y)-u(x))}{|x-y|^{N+s_ip_i}}dy,\,\,\,x\in\mathbb{R}^N,
\end{align}
and the functions $f_i,g_i, i=1,2$ are Carath\'eodory functions, $f_i:\Om \times \R^+ \times \R^+ \to \R^+$, $g_i:\Om \times \R^+ \to \R^+$  satisfying some appropriate conditions which are mentioned in next section. When $s_1=s_2$ the equation reduces to a $(p,q)$ Laplacian problem which appears in a more general reaction-diffusion system 
\begin{equation}\label{intro-1}
    u_t =\text{div}(a(u)\nabla u)+ g(x,u)
\end{equation}
where $a(u)= |\nabla u|^{p-2}\nabla u+ |\nabla u|^{q-2}\nabla u$. Such problems have a wide range of applications in physics and related sciences such as biophysics, plasma physics, and chemical reaction
design, etc. where $u$ describes a concentration, and the first
term on the right-hand side of \eqref{intro-1} corresponds to a diffusion with a diffusion coefficient $a(u)$; the term
$g(x,u)$ stands for the reaction, related to sources and energy-loss processes. A lot of attention has been given to the study of $(p,q)$-Laplace equations in the last few years, for instance refer \cite{deepak,faria-miyagaki, MP, sidi, tanaka, yang-yin}.

Recently, the research community of partial differential equations has been attracted towards the study of fractional $(p, q)$-Laplacian problems and hence, a good amount of literature has been constructed related to this. We start from the article of Chen and Bao \cite{chen-bao} where they studied existence, nonexistence and multiplicity of the following $(p,q)$-fractional Laplacian equation over $\mathbb R^N$,
\begin{align*} &(-\Delta)^{s}_{p}u + a(x)\vert u\vert ^{p-2}u +(-\Delta)^{s}_{q}u + b(x)\vert u\vert ^{q-2}u +\mu(x)\vert u\vert ^{r-2}u \\ &\quad= \lambda h(x)\vert u \vert ^{m-2}u,\quad x\in {\mathbb {R}}^{N}, \end{align*}
with appropriate assumptions on the variables and functions. Next, Bhakta and Mukherjee \cite{Bhakta-Mukh}  studied the following problem in a bounded domain
\begin{align*}
 (-\Delta u)^{s_1}_{p} u +(-\Delta u)^{s_2}_{q} u &= \theta V(x)|u|^{r-2}u  + |u|^{p^*_{s_1}-2}u+ \lambda f(x,u) \;\text{ in }\Om,\\
 u &=0\;\text{ in }\; \mathbb R^N \setminus \Om,
\end{align*}
where $0 < s_2 < s_1 < 1 < r < q < p < N/s_1$, $p^*_{s_1}=\frac{Np}{N-s_1p}$, $s_i\in (0,1)$ for $i=1,2$ and $V$ and $f$ are some appropriate functions. They proved
that there exist weak solutions of the above problem for some range of $\lambda,\theta$. Also, for $V (x) \equiv 1, \lambda = 0$ and
assuming certain other conditions on $n, q, r$, they proved the existence of $cat_{\Om}(\Om)$ non negative solutions
by using Lusternik-Schnirelmann category theory. Using the Nehari manifold technique, Goel et. al \cite{divya1} proved multiplicity results for the following problem over bounded domain
\begin{align}
     (-\Delta u)^{s_1}_{p} u +(-\Delta u)^{s_2}_{q} u &=\lambda a(x)|u|^{\delta-2}u  + b(x)|u|^{r-2}u \;\text{ in }\Om, \label{intro-2}\\
 u &=0\;\text{ in }\; \mathbb R^N \setminus \Om, \nonumber
\end{align}
where $1<\delta \leq q\leq p<r \leq p^*_{s_1}$, $0<s_2<s_1<1$, $N>ps_1$, $\beta>0$ and $a,b$ are sign changing functions. Under appropriate conditions on the parameters, they discuss both sub critical and critical case in the article and also establish that
any weak solution of \eqref{intro-2} belongs to $L^\infty(\Om)\cap C_{loc}^{0,\alpha}(\Om)$ for $\alpha \in (0,1)$ when $2\leq q\leq p<r<p^*_{s_1}$. Regularity results for the equation 
\begin{equation*}
     (-\Delta u)^{s}_{p} u +(-\Delta u)^{s}_{q} u =f(x,u)\;\text{in}\; \mathbb R^N,
\end{equation*}
where  $0<s<1$ and $2\leq q\leq p< N/s$ has been studied in \cite{abreu}. Alves, Ambrosio and Isernia dealt with the following class of problems
\begin{equation*}
    (-\Delta u)^{s}_{p} u +(-\Delta u)^{s}_{q} u +V(\epsilon x)(|u|^{p-2}u+|u|^{q-2}u) =f(u)\;\text{in}\; \mathbb R^N,
\end{equation*}
in \cite{alves}, where $0<s<1$ and $2\leq q\leq p< N/s$. Imposing appropriate assumptions on $V$ and $f$, authors proved existence of ground state solution concentrating on a minimum point of $V$, multiplicity of solutions using Lusternik-Schnirelmann category theory and boundedness of solutions. Very recently, an eigenvalue problem for fractional $(p-q)$ Laplace operator has been studied by Nguyen and Vo in \cite{NV} which is of independent interest. A class of  variable exponent $(p,q)$-fractional Laplacian problems with variable exponents and indeﬁnite weights has been studied in \cite{chung-toan}. We also cite \cite{ambrosio,zhi-yang} as some recent articles in this context for interested readers.

The essence of our article lies in the fact that we study a system of equation with fractional $(p,q)$-Laplace operator over a bounded domain which is first of its kind in literature where we establish existence of three non trivial solutions to $(\mathcal P)$. {The approach is variational and inspired by the articles \cite{bvp,zhang-NA}} and adapted in the fractional framework. In the same article, we motivate that this problem can be extended to considering  more general non local operator $\mathcal L_{\phi_i}$ instead of $(-\Delta)^{s_i}_{p_i}$.

{{{This article has been fragmented into five sections-}}
Section $2$ contains preliminaries, assumptions on $f_i$ and $g_i$ and statements of main results of our article. Proof of Theorem \ref{thm.2} is given in Section 3. We have proved Theorem \ref{Thm.1} in Section 4. Lastly, section $5$ comprises proof of Theorem \ref{Thm.2}.

\noi \textbf{Notations-} $a\vee b= \max\{a,b\}$,\; $a \wedge b=\min \{a,b\}$, \; $|u|_r=\left(\int_\Om |u(x)|^r\,dx\right)^\frac{1}{r}$ for any $r>1$.

\section{Preliminaries}
In this section, we define appropriate function spaces which are required for our analysis. Let $ p> 1,\,s\in(0,1),\, N>ps,\, p_s^*:=\frac{Np}{N-sp}.$ We denote the standard fractional Sobolev space by $W^{s,p}(\Omega)$ endowed with the norm
$$
\|{u}\|_{W^{s,p}(\Om)}:=\|{u}\|_{L^p(\Om)}+\left(\int_{\Om\times\Om} \frac{|u(x)-u(y)|^p}{|x-y|^{N+sp}}dxdy\right)^{1/p}.
$$
We set $Q:=\R^{2N}\setminus (\Om^c \times \Om^c)$, where $\Om^c=\Rn \setminus \Om$ and define $$
X_{s,p}(\Om):=\Big\{u:\mathbb{R}^N\to\mathbb{R}\mbox{ measurable }\Big|u|_{\Omega}\in L^p(\Omega)\mbox{ and }
\int_{Q} \frac{|u(x)-u(y)|^p}{|x-y|^{N+sp}}dxdy<\infty\Big\}.
$$
The space $X_{s,p}(\Om)$ is endowed with the norm defined as
$$\|u\|_{s,p}:=|u|_p+\left(\int_{Q} \frac{|u(x)-u(y)|^p}{|x-y|^{N+sp}}dxdy\right)^{1/p}.$$
We note that in general $W^{s,p}(\Om)$ is not same as $X_{s,p}(\Om)$ as $\Om\times\Om$ is strictly contained in $Q$.
We define the space $X_{0,s,p}(\Om)$ as
$$X_{0,s,p}(\Om) :=\Big\{u \in X_{s,p} : u=0 \quad\text{a.e. in}\quad \Rn \setminus \Om\Big\} $$
or equivalently
as $\overline{C_0^\infty(\Om)}^{X_{s,p}(\Om)}$. It is well-known that for $p>1$,  $X_{0,s,p}(\Om)$ is a uniformly convex Banach space endowed with the norm
$$\|u\|_{0,s,p}=\left(\int_{Q} \frac{|u(x)-u(y)|^p}{|x-y|^{N+sp}}dxdy\right)^{1/p}.$$	

Since $u=0$ in $\Rn\setminus\Om,$ the above integral can be extended to all of $\mathbb{R}^N.$ The embedding
$X_{0,s,p}(\Om)\hookrightarrow L^r(\Om)$ is continuous for any $r\in[1,p^*_s]$ and compact for $r\in[1,p^*_s).$
{ Moreover, for $1<q\leq p$, $X_{0,s_1,p}(\Om)\subset X_{0,s_2,q}(\Om)$ (see Lemma 2.2 in Section 2 of \cite{Bhakta-Mukh}). 

Let us recall some topological tools which will be used to study our problem $(\mathcal{P})$.
\begin{definition}
(Nehari Manifold):
	Let $J \in C^1(X,\R)$ be such that $J'(0)=0$, then the constraint set
	$$ \mathcal{N}:=\{ u \in X: \<J'(u),u\>=0,\,u\neq 0   \}
	$$
	is called a Nehari manifold of $X$.
	\end{definition}


\begin{definition}
Let $X$ be a Banach space, $\eta:X\to \R$ be a function. For $c\in \R,$ let us consider the punctured level set of $\eta$ at $c$ by $$L_\eta^c = \big\{u \in X : \eta(u)=c,u \neq 0 \big\}. $$ We say $L_\eta^c$
	has the sphere property, if the following hypotheses are satisfied:
	\begin{itemize}
	\item [(i)]
		$\eta\in C(X);$
		\item [(ii)]
		there is a homeomorphic mapping between $L_\eta^c$ and the unit sphere of $X$;
		\item [(iii)]
		for any fixed $u \in X \setminus \{0\}$, there exists a unique $t_u\in (0,\infty)$ such that $f(t_uu)=c;$
		\item [(iv)]
		$X$ is separated into two open connected subsets by $L_\eta^c$ and the origin is contained in one of these subsets.
	\end{itemize}	
\end{definition}	
\begin{theorem}\label{thm.link}
	Let $X,Y$ be Banach spaces with the following direct sum decomposition:
	$$
	X=X_1 \oplus X_2, Y=Y_1 \oplus Y_2,
	$$
	where $X_1,Y_1$ are finite dimensional subspaces of $X,Y$ respectively. Let $\eta\in C(X),\kappa \in C(Y),\,\,c,d\in\R$ and $L_\eta^c, L_\kappa^d$ have the sphere property. Let $(e_x,e_y) \in X_2 \times Y_2$ such that $\eta(e_x)>c$ and $\kappa(e_y)>d$. Let us denote
	\begin{gather*}
		\mathcal{Q}^X=\big\{ u+te_x: u \in X_1 \cap B_{R_1}, \, t \in [0,1]   \big\},
	\end{gather*}	
	\begin{gather*}
	\mathcal{Q}^Y=\big\{ v+s e_y: v \in Y_1 \cap B_{R_2}, \, s \in [0,1]   \big\},
\end{gather*}	
$\mathcal{N}_1^{X_2}=L_\eta^c \cap X_2$, $\mathcal{N}_2^{Y_2}=L_\kappa^d \cap Y_2$,
$\mathcal{Q}=\mathcal{Q}^X \times \mathcal{Q}^Y$ and $\mathcal{N}=\mathcal{N}_1^{X_1} \times \mathcal{N}_2^{X_2}$. Then, $\pa Q$ links $\mathcal{N}$.
	\end{theorem}	

\subsection{Assumptions on $f_i$ and $g_i$, $i=1,2$}
 We assume the following:
\begin{itemize}
\item [\bf(A1)]
	There exists a function $F \in C^1(\Om \times \R \times \R,\R )$ such that
	\begin{gather*}
	 \left(  \frac{\pa F}{\pa u}(x,u,v),  \frac{\pa F}{\pa v}(x,u,v)  \right)=\big(  f_1(x,u,v), f_2(x,u,v)  \big),\text{ for all } (x,u,v) \in \Om \times \R \times \R.
	\end{gather*}	
\item [\bf(A2)]
For $i=1,2$,\; \,$g_i \in C(\bar{\Om} \times \R^+,\R^+)$ and there exist constants
$q_i \in (p_i,p_i^*)$ and $C_i>0$ such that
	    $$ |g(x,u)| \leq C_i \big( |u|^{p_i-1}+|u|^{q_i-1}  \big) \, \text{for all}\, (x,u) \in \Om \times \R^+.
	$$
\item [\bf(A3)]	
There exists constant $\al_i>p_1\vee p_2$ such that
\begin{gather*}
0<\al_i G_i(x,u) \leq u g_i(x,u)\,\text{for}\, (x,u)\in \Om \times (\R^+ \setminus \{0\}),	
\end{gather*}		
where $G_i(x,u)=\int_0^u g_i(x,\tau)\,d\tau;i=1,2$.
\item[\bf(A4)]
There holds for $i=1,2$,
\begin{gather*}
	\limsup_{u \to 0}\frac{g_i(x,u)}{|u|^{p_i-2}u}<\la_{1,p_i} \,\text{uniformly w.r.t}\;\, x \in \Om,
\end{gather*}	
where $\la_{1,p_i}$ are the first eigenvalue of $(-\De)^{s_i}_{p_i}$ in $X_{0,s_i,p_i}$, that is,
$$
\la_{1,p_i}=\inf_{u \in X_{0,s_i,p_i} \setminus \{0\}} \frac{\|u\|_{0,s_i,p_i}^{p_i}}{|u|_{p_i}^{p_i}}.
$$
\item[\bf(A5)]There holds,
\begin{gather*}
	\frac{g_i(x,u)}{u_i^{p-1}} \,\text{is an increasing function of}\, u \in \R^+ \setminus \{0\}.
\end{gather*}	
\item[\bf(A6)]
There exists constant {$q \in (1,p_1^* \wedge p_2^*)$} and $C_3>0$ such that
\begin{equation*}
	|f_1(x,u,v)|+|f_2(x,u,v)| \leq C_3 \big( |u|^{q-1}+|v|^{q-1}+1 \big).
\end{equation*}	
\item[\bf(A7)]
There exist $\ba_i \in (1, \al_1 \wedge \al_2)$ for $i=1,2$, $C_4>0,R>0$ such that,
\begin{equation*}
	u f_1(x,u,v)+v f_2(x,u,v) \leq C_4  \big( |u|^{\ba_1}+|v|^{\ba_2} \big)\,\text{for}\, |u|+|v|\geq R.
\end{equation*}	

\item[\bf(A8)]
	There holds 
	\begin{gather*}
		\limsup_{|u|_{p_1}^{p_1}+|v|_{p_2}^{p_2} \to \infty} \frac{(p_1\vee p_2)F(x,u,v)}{|u|^{p_1} +|v|^{p_2}} \leq h(x)\,\,\text{uniformly for a.e.}\, x \in \Om.
	\end{gather*}
for some $h \in L^{\infty}(\Om)$.
\end{itemize}	
{ A model example satisfying (A1) to (A8) can be taken as 
 $$F(x,u,v)=u^{q-2}v^2+v^{q-2}u^2,$$
for $x\in \Om$ and $u,v>0$. Then, we have,
 $$f_1(x,u,v)=\frac{\pa F}{\pa u}(x,u,v)=(q-2)u^{q-3}v^2+2u v^{q-2},$$ and
 $$f_2(x,u,v)=\frac{\pa F}{\pa v}(x,u,v)=2v u^{q-2}+(q-2)v^{q-3}u^2. $$
  Choosing $q_i \in (p_1 \vee p_2, p_1^*\wedge p_2^*)$ for $i=1,2$, we can take $g_1(x,u)= u^{q_1-1}$ and $g_2(x,v)=v^{q_2-1}$ for $x\in \Om$ and $u,v>0$.}
Our first main result in the article is the following:

\begin{theorem}\label{thm.2}
	Let $f_i,g_i$ satisfy (A1)-(A7) for $i=1,2$, together with the condition
	\begin{gather}\label{eq:*}
		f_i(x,u,0)=f_i(x,0,v)=0,\,\text{for}\,\,i=1,2,
		\end{gather}
	for a.e. $x \in \Om$, for all $u,v \in \R^+$. Also, we assume $F \in C^1(\Om \times \R^+ \times \R^+,\R )$ in (A1). Then, there exists $\La>0$ such that for any $\la \in (0,\La)$, system $(\mathcal{P})$ has atleast three non-negative solutions. Moreover, if both the problem, for $i=1,2$,
	 \begin{equation*}
	 	(\mathcal{P'})
	 	\left\{\begin{aligned}
	 		(-\Delta)^{s_i}_{p_i} u  &=g_i(x,u)  \,\mbox{ in }\, \Om, \\
	 		u&= 0  \,\mbox{in }\, \Rn \setminus \Om, \\
	 	\end{aligned}
	 	\right.
	 \end{equation*}
 have a unique positive solution, then for any $\la \in (0,\La)$, system $(\mathcal{P})$ has at least three non-negative solutions; among them, one is nontrivial positive solution.
\end{theorem}	

We shall prove the existence of weak solutions of $(\mathcal{P})$ for $\lambda=1, \,\, g_1\equiv0\equiv g_2$ by means of variational methods in the space $X := X_{0,s_1,p_1}(\Om) \times X_{0,s_2,p_2}(\Om)$ endowed with the norm given by $\|(u,v)\|= \|u\|_{0,s_1,p_1}+ \|v\|_{0,s_2,p_2}$ for all $(u,v)\in X$.
{{
\begin{definition}
 An element $(u,v) \in X$ is said to be a weak solution of $(\mathcal{P})$
 if $(u,v)$ satisfies
 $$ \<J'(u,v),(w,z)\>=0 \, \mbox{for all}\, (w,z) \in X.
 $$
\end{definition}
Our second main out-turn in the article is the following.

\begin{theorem}\label{Thm.1} Let $\lambda=1, \,\, g_1\equiv0\equiv g_2$
and $f_1,f_2$ satisfy (A1) with assumption $F \in C^1(\Om \times \R^+ \times \R^+,\R )$, (A6), (A8) and $h(x) \leq \bar{\la}$ in $\Om$ and $h(x) < \bar{\la}$ on a subset of $\Om$ with positive measure, $\bar{\la}=min\{\la_{1,p_1},\la_{1,p_2} \}$; $\la_{1,p_i}$'s are defined in (A4). Then, there exists a weak solution $(u,v) \in (X_{0,s_1,p_1}(\Om) \times X_{0,s_2,p_2}(\Om))$ of $(\mathcal{P})$.
\end{theorem}

}}

\begin{remark}
 We may consider more general nonlocal operator, for example,
 we consider the following system of non-local quasilinear elliptic equations:
\begin{equation*}
	(\mathcal{Q})
	\left\{\begin{aligned}
		-\mathcal{L}_{\phi_1} u  &=f_1\,(x,u,v) \,\mbox{ in }\, \Om, \\
		-\mathcal{L}_{\phi_2} v &=f_2\,(x,u,v) \,\mbox{ in }\, \Om, \\
		u=v &= 0  \,\mbox{in }\, \Rn \setminus \Om, \\
	\end{aligned}
	\right.
\end{equation*}
where the operators $\mathcal{L}_{\phi_i}, i=1,2$ are defined by:
\begin{equation*}
	\<-\mathcal{L}_{\phi_1}u,w\>=\int_{\R^{2N}} \phi_1 (u(x)-u(y))(w(x)-w(y))K_1(x,y)\,dxdy,
\end{equation*}
and
\begin{equation*}
	\<-\mathcal{L}_{\phi_2}u,w\>=\int_{\R^{2N}} \phi_2 (u(x)-u(y))(z(x)-z(y))K_2(x,y)\,dxdy,
\end{equation*}
for all $w,z \in C_c^{\infty}(\Om)$, the functions $\phi_i,\; i=1,2$ are assumed to be continuously differentiable satisfying $\phi_i(0)=0,i=1,2$,
\begin{gather}\label{*1}
  \text{the function }\,  t \mapsto t \phi_i(t) \, \text{ is convex},
\end{gather}
and there exists $\gamma_i>0,i=1,2$ such that
\begin{gather*}
	\frac{1}{\gamma_i}|t|^{p_i} \leq \phi_i(t)t \leq \gamma_i|t|^{p_i}, \, \text{for all }\, t \in \R,
\end{gather*}
and $K_i:\Rn \to \R$ are assumed to be measurable, symmetric and satisfy for some $\delta_i\geq 1,$
\begin{gather}\label{*2}
	\frac{1}{\delta_i|x-y|^{N+s_ip_i}} \leq K_i(x,y) \leq \frac{\delta_i}{|x-y|^{N+s_ip_i}},
\end{gather}
for all $x,y \in \Rn$.
With this in hand, we have the following outcome.
\begin{theorem}\label{Thm.2}Let $\lambda=1, \,\, g_1\equiv0\equiv g_2$
and $f_1,f_2$ satisfy (A1),(A6), (A8) and $h(x) \leq \frac{\bar{\la}}{\max\{\gamma_1,\gamma_2\}}$ in $\Om$ and $h(x) < \frac{\bar{\la}}{\max\{\gamma_1,\gamma_2\}}$ on a subset of $\Om$ with positive measure, $\bar{\la}=min\{\la_{1,p_1}^\mathcal{L},\la_{1,p_2}^\mathcal{L} \}$; $\la_{1,p_i}^\mathcal{L}$'s are defined by
$$
\la_{1,p_i}^\mathcal{L}=\inf_{u \in X_{0,s_i,p_i} \setminus \{0\}} \frac{\langle -\mathcal{L}_{\phi_i}u,u\rangle}{{|u|_{p_i}}^{p_i}}.
$$
Then, there exists a weak solution $(u,v) \in (X_{0,s_1,p_1}(\Om) \times X_{0,s_2,p_2}(\Om))$ of $(\mathcal{Q}).$
\end{theorem}
\end{remark}

\section{Proof of Theorem \ref {thm.2}}
\noindent For $i=1,2,$ let us first define the functional $J_i: X_{0,s_i,p_i} \to \R$ by
\begin{equation*}
	J_i(u)=\frac{1}{p_i}\|u\|_{0,s_i,p_i}^{p_i}-\Iom G_i(x,u)\, dx,\,\, u\in X_{0,s_i,p_i}
\end{equation*}	
where $G_i$'s are defined in (\textit{A3}).
Let us denote the Nehari manifold of $J_i$ on $X_{0,s_i,p_i}$ by $\mathcal{N}_i$ for $i=1,2$ that is
$$\mathcal N_i = \left\{u \in X_{0,s_i,p_i}\setminus\{0\}:\; \|u\|_{0,s_i,p_i}= \int_\Om g_i(x,u)u~dx\right\}.$$
We prove the following  result.
\begin{theorem}\label{thm.1}
	Let $g_i$ for $i=1,2$ satisfy (A2)-(A5). For $i=1,2$, let $e_i \in X_{0,s_i,p_i}$
	be such that {$J_i(e_i)> 0$}. Let us define
$$	\mathcal{Q}:=\big\{ \la(e_1,0)+(1-\la)(0,e_2) : \la \in [0,1] \big\} \subset X$$
and $\mathcal{N}=\mathcal{N}_1 \times \mathcal{N}_2$, where $\mathcal{N}_i$'s are defined above. Then $\pa \mathcal{Q}$ links $\mathcal{N}$.
\end{theorem}	

\begin{proof}
    We prove this result in two steps.
	
\textbf{Step-1}. The Nehari manifold $\mathcal{N}_i$'s both have the sphere property.	
 {{
 The proof is similar to the proof of Lemma 4.1 of \cite{Willem}, on page 72, where the functional $\varphi$ is replaced with $J_i.$ So we omit it.
}}

\allowdisplaybreaks

 \textbf{Step-2}. We note that $J_i$'s satisfy the following for $i=1,2$,
\begin{align*}
J_i(u)=\begin{cases}
	0,\,\,u=0,\\
	>0, \,u=e_i.
	\end{cases}
\end{align*}	
By assumption {(A3), we can get that $u^{\alpha_i} \leq C G_i(x,u)$ for some constant $C>0$ where $(x,u)\in \Om \times \mathbb R^+$. This gives that $\lim\limits_{t\to \infty}J_i(te_i)=-\infty$ since $\alpha_i\in \max\{p_1,p_2\}$ from which} we note that there exists $t_i \in (0,1)$ such that
\begin{align*}	
	\<J_i'(te_i),e_i\>=\begin{cases}
	0,\,t=t_i,\\
	>0, 0<t<t_i,\\
	<0,\,t_i<t<\infty.
	\end{cases}
\end{align*}	
Hence, we have, $\<J'_i(e_i),e_i\><0$. Using Step-1 along with Theorem \ref{thm.link}
we conclude that $\pa Q$ links $\mathcal{N}$.
\end{proof}

With this machinery in hand, we commence to prove our first main result.
\begin{proof}[\textbf{Proof of Theorem \ref{thm.2}}] 
{{
We prove this result in three steps. To prove this result, we will use Theorem \ref{thm.1} to the functional $J_i$.
In the first step, we show that $J$ satisfies Palais-Smale condition. Then, in the next step, for the functional $J_i$, we consider the few critical levels and values corresponding to $J_i$. With these critical values in hand, we prove that infimum over such Nehari manifolds for the functionals $J_1$ and $J_2$ are achieved and they form two solutions to our problem. Using these, we will construct a set $\mathcal{Q}$ (defined in Theorem \ref{thm.1}) such that $\mathcal{Q}$ links $\mathcal{N}_i$. Lastly, we prove the existence of our third solution using our step-(1). This yields our result.
}}

	\textbf{Step-1}. Let us consider the cone $X^+= \{(u,v)\in X:\; u,v\geq 0\}$ and  define the following $J:X^+ \to \R$ by
	\begin{align}\label{eq:J}
		J(u,v)=J_1(u)+J_2(v)-\la\Iom F(x,u,v)\,dx \;\;\,\text{for all}\, (u,v) \in X^+.
	\end{align}	
	In this step, we show that $J$ satisfies (PS) condition.
	Let $\{(u_n,v_n)\} \subset X^+$ be a sequence such that
	\begin{align}\label{eq:PS}
		\{J(u_n,v_n) \}_{n \geq 1} \, \text{is bounded},\, J'(u_n,v_n) \to 0,\, n \to \infty.
	\end{align}
	Then there exists $C_5>0$ such that $J(u_n,v_n) \leq C_5\; \, \text{for all}\, n \in \N$, that is,
	\begin{align}\label{eq:bound}
		J_1(u_n)+J_2(v_n) -\la \Iom F(x,u_n,v_n)\,dx \leq C_5 \, \text{for all} \, n\in \N,
	\end{align}	
	and
	\begin{gather*}
	\bigg|  \frac{\<J'(u_n,v_n),(u_n,v_n)\>}{\| (u_n,v_n)\|}  \bigg| <1 \,\text{for large}\, n \in \N,
	\end{gather*}	
	that is,
	\begin{align}\label{eq:J'}
		|\<J'(u_n,v_n),(u_n,v_n)\>| \leq \|(u_n,v_n)\| \,\text{for large}\, n \in \N.
	\end{align}
{We first realise that using $(A7)$ and continuity of $f_i$, we get
	\begin{equation}\label{sys-1}
	\begin{split}
	    F(x,u,v) &= F(x,0,0) + \int_0^1 \frac{d}{dt} F(x,tu,tv)~dt\\
	    &= F(x,0,0) + \left(\int_0^{\frac{R}{|u|+|v|}}+ \int^{1}_{\frac{R}{|u|+|v|}}\right)(f_1(x,tu,tv)u+ f_2(x,tu,tv)v)~dt\\
	    & \leq F(x,0,0)+ C^\prime + C_4(|u|^{\beta_1}+ |v|^{\beta_2}),\; \text{for}\; (x,u,v)\in  \Om \times \mathbb{R}^+ \times \mathbb{R}^+ 
	    \end{split}
	\end{equation}
	for some constant $C^\prime>0$.}
	Let $r \in ({p_1 \vee p_2}, \al_1 \wedge \al_2)$. Then using the above estimate  with (\ref{eq:PS}), (\ref{eq:bound}) and the assumptions (A1)-(A7), we obtain for large $n$,
	\begin{equation*}
		\begin{aligned}
			C_5+\frac{1}{r}\|(u_n,v_n)\| &\geq J(u_n,v_n)-\frac{1}{r}\<J'(u_n,v_n),(u_n,v_n)\>\\
			&=J_1(u_n)+J_2(u_n) -\la \Iom F(x,u_n,v_n)\,dx\\
			&\quad-\frac{1}{r}\bigg\{ \|u_n\|_{0,s_1,p_1}^{p_1}+\|v_n\|_{0,s_2,p_2}^{p_2}-\la \Iom \big[ u_n f_1(x,u_n,v_n)+v_n f_2(x,u_n,v_n)  \big] \,dx\\
			& \quad-\Iom \big[ u_n g_1(x,u_n)+v_n g_2(x,v_n)  \big]\,dx  \bigg\}\\
			&=\left( \frac{1}{p_1}-\frac{1}{r} \right) \|u_n\|_{0,s_1,p_1}^{p_1}+\left( \frac{1}{p_2}-\frac{1}{r} \right) \|v_n\|_{0,s_2,p_2}^{p_2}\\
			&\quad+\la \bigg[
			\Iom \frac{1}{r} \big( u_n f_1(x,u_n,v_n)+v_n f_2(x,u_n,v_n)\big) \,dx -\Iom F(x,u_n,v_n)\,dx \bigg]\\
			&\quad-\bigg[ \Iom \big( G_1(x,u_n)-\frac{1}{r}\, u_n \, g_1(x,u_n)  \big)\,dx
			+\Iom \big( G_2(x,v_n)-\frac{1}{r}\,v_n\, g_2(x,v_n)  \big)\,dx  \bigg]\\
			&\geq \left(\frac{1}{p_1}-\frac{1}{r}\right) \|u_n\|_{0,s_1,p_1}^{p_1}
			+\left(\frac{1}{p_2}-\frac{1}{r}\right)\|v_n\|_{0,s_2,p_2}^{p_2} \\
			& \quad+\la \bigg[
			\Iom \frac{1}{r} \big[ u_n f_1(x,u_n,v_n)+v_n f_2(x,u_n,v_n)\big] \,dx -\Iom F(x,u_n,v_n)\,dx \bigg]\\
			& \quad {+ \frac{(\al_1-r)}{r}\Iom G_1(x,u_n)\,dx+\frac{(\al_2-r)}{r}\Iom G_2(x,v_n)\,dx}\\
			&\geq \left(\frac{1}{p_1}-\frac{1}{r}\right) \|u_n\|_{0,s_1,p_1}^{p_1}
			+\left(\frac{1}{p_2}-\frac{1}{r}\right)\|v_n\|_{0,s_2,p_2}^{p_2} +
			C_6 ( |u_n|_{\al_1}^{\al_1}+|v_n|_{\al_2}^{\al_2} )\\
			&\quad -\la C_7( |u_n|_{\ba_1}^{\ba_1}+|v_n|_{\ba_2}^{\ba_2} )-C_8\\
		\end{aligned}
	\end{equation*}
	This implies,
	\begin{equation}\label{eq:u_n,v_n}
		C_5+\frac{1}{r}\|(u_n,v_n)\| \geq \left(\frac{1}{p_1}-\frac{1}{r}\right) \|u_n\|_{0,s_1,p_1}^{p_1}
		+\left(\frac{1}{p_2}-\frac{1}{r}\right)\|v_n\|_{0,s_2,p_2}^{p_2}-C_{9}.
	\end{equation}
{It is easy to verify that $a^{p_1}+b^{p_2} \leq 2(a+b)^{\max\{p_1,p_2\}}$ for any $a,b\in \mathbb{R}^+$ which applied to (\ref{eq:u_n,v_n}) yields that $\{(u_n,v_n)\}$ is bounded in $X^+$ since $\min\{p_1,p_2\}>1$}. Therefore, up to a subsequence, we may assume that there exists $(u,v) \in X^+$(since $X^+$ is a closed subspace of $X$) such that
	$(u_n,v_n) \deb (u,v)$ weakly in $X^+$, $u_n \to u$ strongly in $L^{\ga_1}(\Rn)$,  $v_n \to v$ strongly in $L^{\ga_2}(\Rn)$ for $\ga_i \in [1,p_i^*)$, $i=1,2$. Also we assume that $(u_n,v_n) \to (u,v)$ as $n\to \infty$ pointwise a.e. in $\Om$.
	This immediately implies,
	\begin{gather*}
		\left| \Iom g_1(x,u_n) (u_n-u)\,dx\right| \leq |g_1(\cdot, u_n)|_{q_1'}|u_n-u|_{q_1}
		\to 0 \, \text{as}\, n \to \infty,
	\end{gather*}	
	where $q_1'=\frac{q_1}{q_1-1}$.
We also note that,	
\begin{gather*}
	\left|\Iom g_2(x,v_n) (v_n-v)\,dx \right| \leq |g_2(\cdot,v_n)|_{q_2'}|v_n-v|_{q_2} \to 0 \, \text{as}\, n \to \infty,
\end{gather*}		
	\begin{gather*}
		\left|\Iom f_1(x, u_n, v_n) (u_n-u)\,dx \right| \leq |f_1(\cdot, u_n, v_n)|_{q}|u_n-u|_{q} \to 0 \, \text{as}\, n \to \infty,
	\end{gather*}	
	\begin{gather*}
	\left|\Iom f_2(x, u_n, v_n) (v_n-v)\,dx \right| \leq |f_2(\cdot, u_n, v_n)|_{q}|v_n-v|_{q} \to 0 \, \text{as}\, n \to \infty.
\end{gather*}	
These together with (\ref{eq:PS}) implies,
\begin{equation*}
\begin{aligned}
&\Iom \frac{|u_n(x)-u_n(y)|^{p_1}}{|x-y|^{N+s_1 p_1}}\,dx dy
+\Iom \frac{|v_n(x)-v_n(y)|^{p_2}}{|x-y|^{N+s_2 p_2}}\,dx dy\\
&-\Iom \frac{|u_n(x)-u_n(y)|^{p_1-2} ( u_n(x)-u_n(y))( u(x)-u(y) )}{|x-y|^{N+s_1 p_1}}\,dx dy\\
&-\Iom \frac{|v_n(x)-v_n(y)|^{p_2-2} ( v_n(x)-v_n(y))( v(x)-v(y) )}{|x-y|^{N+s_2 p_2}}\,dx dy\\
& \to 0\,\text{as}\, n \to \infty,
\end{aligned}
\end{equation*}
that is,
\begin{gather}\label{eq:-i}
	\< (u_n,v_n), (u_n,v_n)-(u,v)\>_{X,X'} \to 0\,\text{as}\, n \to \infty.
\end{gather}	
Hence, using (\ref{eq:-i}), we obtain,
\begin{equation*}
\| (u_n,v_n)-(u,v) \| \to 0 \,\text{as}\, n \to \infty.
\end{equation*}	
which implies that $(u_n,v_n) \to (u,v)$ strongly in $X$ as $n \to \infty$. Therefore, we conclude that $J$ satisfies the (PS) condition.

\textbf{Step-2}. Let us define for $i=1,2$,
\begin{gather*}
	c_i^*:=\inf_{u \in \mathcal{N}_i} J_i(u), \, 	c_i^{**}:=\inf_{u \in  X_{0,s_i,p_i}\setminus \{0\}  } \max_{t \geq 0} J_i(tu), \,
	c_i=\inf_{\ga \in \Ga_i} \max_{t \in [0,1]} J_i(\ga(t)) \,\text{where}\,
\end{gather*}	
$$
\Ga_i:=\bigg\{ \ga \in C([0,1] , X_{0,s_i,p_i}) : \ga(0)=0, J_i(\ga(1))<0  \bigg\}.
$$
 We have the following: $c_i$ is a critical value of $J_i$ and $c_i^*=c_i^{**}=c_i>0$. The proof is similar to the proof of Theorem 4.2 of \cite{Willem}, on page 73. So we omit it.


\textbf{Step-3}. From Step-2, we have that $c_i$'s are the critical values of $J_i$
for $i=1,2$. Therefore, there exists $\bar{u} \in \mathcal{N}_1$ and $\bar{v} \in \mathcal{N}_2$ such that
\begin{gather*}
	J_1(\bar{u})=c_1,\, J_2(\bar{v})=c_2 \,\text{and}\, J_1'(\bar{u})=0=J_2'(\bar{v}).
\end{gather*}	
Hence, we obtain, $(\bar{u},0)$ and $(0,\bar{v})$ as {non-negative solutions} of
 \begin{equation}\label{eq:P}
	\left\{\begin{aligned}
		(-\Delta)^{s_1}_{p_1} \bar{u}  &=g_1(x, \bar{u})  \,\mbox{ in }\, \Om, \\
		(-\Delta)^{s_2}_{p_2} \bar{v}  &=g_2\,(x,\bar{v}) \,\mbox{ in }\, \Om, \\
		\bar{u}= 0&=\bar{v}  \,\mbox{in }\, \Rn \setminus \Om, \\
	\end{aligned}
	\right.
\end{equation}
Since (\ref{eq:*}) holds, we conclude that $(\bar{u},0)$ and $(0,\bar{v})$ are non-negative solutions of  $(\mathcal{P})$ with
\begin{gather*}
	J(\bar{u},0)=J_1(\bar{u})=c_1\,\text{and}\,J(0,\bar{v})=J_2(\bar{v})=c_2.
\end{gather*}	
To obtain the third non-negative solution, let us consider the term $J(t\bar{u},s\bar{v})$ for $t,s \geq 1$ given by
\begin{equation*}
	\begin{aligned}
		J(t\bar{u},s\bar{v})&=J_1(t\bar{u})+J_2(s\bar{v})-\la\Iom F(x,t\bar{u},s\bar{v})\,dx\\
		&=\frac{1}{p_1}\|t\bar{u}\|_{0,s_1,p_1}^{p_1}+\frac{1}{p_2}\|s\bar{v}\|_{0,s_2,p_2}^{p_2}-\Iom G_1(x,t\bar{u})\,dx
		-\Iom G_2(x,s\bar{v})\, dx\\
		&\quad-\la \Iom F(x,t\bar{u},s\bar{v})\,dx\\
		&=\frac{t^{p_1}}{p_1}\|\bar{u}\|_{0,s_1,p_1}^{p_1}+\frac{s^{p_2}}{p_2}\|\bar{v}\|_{0,s_2,p_2}^{p_2}-\Iom G_1(x,t\bar{u})\,dx-\Iom G_2(x,s\bar{v})\,dx\\
		&\quad-\la \Iom \int_0^1 \big[  t\bar{u} f_1(x,r t \bar{u}, rs\bar{v}) +s\bar{v}f_2(x,rt\bar{u},rs\bar{v})  \big]\,drdx\\
		&\leq \frac{t^{p_1}}{p_1}\|\bar{u}\|_{0,s_1,p_1}^{p_1}+\frac{s^{p_2}}{p_2}\|\bar{v}\|_{0,s_2,p_2}^{p_2}-C \big( |t\bar{u}|_{\al_1}^{\al_1}+|s\bar{v}|_{\al_2}^{\al_2} \big)+C\\
		& {\leq C (t+s)^{p_1 \vee p_2}-C(t+s)^{\al_1 \wedge \al_2}+C,}
	\end{aligned}	
\end{equation*}	
with some generic constant $C>0$ (independent of $t,s$ but depends on $\bar u$, $\bar v$). Therefore, there exists a $L>0$ such that,
\begin{gather}\label{eq:L}
	J(t\bar{u},s\bar{v})<0 \;\,\text{for all}\;\, s,t \;\,\text{with}\;\, {s+t\geq L}.
\end{gather}	
Let $\mathcal{Q}:=\big\{ (tL\bar{u},sL \bar{v}) \in X: (t,s) \in [0,1] \times [0,1]     \big\}$
and {$\mathcal N=\mathcal{N}_1 \times \mathcal{N}_2$}.
We observe that,
\begin{gather*}
	J_1(L\bar{u})=J(L\bar{u},0)<0 \,\text{and}\,	J_2(L\bar{v})=J(0,L\bar{v})<0.
\end{gather*}	
Hence, applying Theorem \ref{thm.link} we assert that $\pa \mathcal{Q}$ links $\mathcal{N}$.\\
\\
\textbf{Claim}. For $\la>0$ small enough,
\begin{gather*}
	\sup_{\pa \mathcal{Q}} J(u,v) < \inf_{\mathcal{N}} J(u,v).
\end{gather*}	
\textit{Proof of Claim}. We know from the definition of $c_i,c_i^{**}$, as in Step-2,
that
\begin{gather*}
	{\inf_{u \in X_{0,s_i,p_i} \setminus \{0\}} \max_{t \geq 0} J_i(tu)=c_i^{**}=c_i = J_i(\bar u)\leq \max_{t\in [0,1]} J_i(tL\bar u).}
\end{gather*}	
By virtue of (\ref{eq:L}) and since $c_i>0$ from Step-2, we obtain
\begin{gather}\label{eq:J-*}
	\sup_{(u,v) \in \pa \mathcal{Q}} J(u,v) \leq c_1 \vee c_2 \,\text{for all}\, \la>0.
\end{gather}	
{For $(u,v)\in \mathcal N$, using $(A3)$ we see that
$$J_1(u)\geq \left(\frac{1}{p_1}-\frac{1}{\alpha_1}\right)\|u\|^{p_1}_{0,s_1,p_1},\quad J_2(v)\geq \left(\frac{1}{p_2}-\frac{1}{\alpha_2}\right)\|v\|^{p_2}_{0,s_2,p_2} .$$
Using \eqref{sys-1}, we infer that for $(u,v)\in \mathcal N$
\begin{equation*}
    \begin{split}
        \int_\Om F(x,u,v) ~dx &\leq C_{10} + C_4 \left(|u|^{\beta_1}_{\alpha_1}|\Om|^{1-\frac{\beta_1}{\alpha_1}} + |v|^{\beta_2}_{\alpha_2}|\Om|^{1-\frac{\beta_2}{\alpha_2}}\right)\\
        & \leq C_{10} + \left(C_{11}|u|_{\alpha_1}+ C_{12}|v|_{\alpha_2}\right), \;\text{using Young's inequality}\\
        & \leq C_{10}+ \int_\Om G_1(x,u)~dx + \int_{\Om} G_2(x,v)~dx \\
        & \leq C_{10}+ \frac{1}{\alpha_1}\int_{\Om}ug_1(x,u)~dx + \frac{1}{\alpha_2}\int_{\Om}vg_2(x,v)~dx, \;\text{from (A3)}\\
        & = C_{10} + \frac{\|u\|_{0,s_1,p_1}^{p_1}}{\alpha_1}+ \frac{\|v\|_{0,s_2,p_2}^{p_2}}{\alpha_2}
    \end{split}
\end{equation*}
From the above estimates, if we assume $\la <\min \left\{ \frac{\alpha_1-p_1}{p_1}, \frac{\alpha_2-p_2}{p_2}\right\}$ then  for any $(u,v) \in \mathcal{N}$ we have that
\begin{align*}
	J(u,v) &\geq J_1(u) + J_2(u) -\la C_{10} + \frac{\la \|u\|_{0,s_1,p_1}^{p_1}}{\alpha_1}+ \frac{\la \|v\|_{0,s_2,p_2}^{p_2}}{\alpha_2}\\
	&\geq \bigg(  1-\frac{\la p_1 }{\al_1-p_1}  \bigg) J_1(u) + \bigg(  1-\frac{\la p_2 }{\al_2-p_2}  \bigg) J_2(v) -\la \bar{C}\\
	& \geq \bigg(  1-\frac{\la p_1 }{\al_1-p_1}  \bigg) c_1 + \bigg(  1-\frac{\la p_2 }{\al_2-p_2}  \bigg) c_2 -\la \bar{C}
\end{align*}	
for some $\bar{C}>0$. Let us take
\begin{equation*}
\La=\min \bigg\{ \frac{\al_1-p_1}{p_1 },  \frac{\al_2-p_2}{p_2 }, \frac{c_1+c_2-(c_1 \vee c_2)}{\frac{p_1  c_1}{\al_1-p_1} +\frac{p_2  c_2}{\al_2-p_2}+ \bar{C}}   \bigg\}.
\end{equation*}
then it is easy to see that whenever
 $\la \in (0,\La)$,
\begin{gather}
	c_1 \vee c_2 <  \left(  1-\frac{\la p_1 }{\al_1-p_1}  \right) J_1(u) + \left( 1-\frac{\la p_2 }{\al_2-p_2}  \right) J_2(v) -\la \bar{C}\leq J(u,v), 
\end{gather}	}
for all $(u,v)\in \mathcal N$ which implies,
\begin{gather}\label{eq:J-N}
	\inf_{(u,v) \in \mathcal{N}} J(u,v) \geq c_1 \vee c_2.
\end{gather}	
Using (\ref{eq:J-N}), we have from (\ref{eq:J-*}) that,
$$
\sup_{(u,v) \in\; \pa \mathcal{Q}} J(u,v) \leq  \inf_{(u,v) \in \mathcal N } J(u,v).$$
Let us define $\bar{c}:=\inf\limits_{\ga \in \Ga} \sup\limits_{(u,v) \in \mathcal{Q}} J(\ga(u,v))$,
where
$$\Ga:=\bigg\{  \ga \in C(\mathcal{Q}, X): \ga|_{\pa \mathcal{Q}}=Id|_{\pa \mathcal{Q}}   \bigg\}.$$
\begin{equation*}
	\bar{c} \geq \inf_{(u,v) \in \mathcal{\tilde{N}}} J(u,v) > c_1 \vee c_2
\end{equation*}	
and $\bar{c}$ is a critical point of $J$. Hence,
{{$(\mathcal{P'})$}} has at least three non-negative solutions for $\la \in (0,\La)$. Furthermore, if the system {{$(\mathcal{P'})$}} have a unique positive solution for $i=1,2$, the third non-negative solution of {{$(\mathcal{P'})$}} is positive because if one of the component of third solution is zero, then the value of $J$ is either $c_1$ or $c_2$ which is a contradiction.
\end{proof}	
\section{Proof of Theorem \ref{Thm.1}}
\noi This section consists of the study of $(\mathcal P)$ with $g_1\equiv 0 \equiv g_2$ and $\la =1$ that is
\begin{equation*}
	(\mathcal{P}_0)
	\left\{\begin{aligned}
		(-\Delta)^{s_1}_{p_1} u  &=f_1\,(x,u,v) \,\mbox{ in }\, \Om, \\
		(-\Delta)^{s_2}_{p_2} v  &=f_2\,(x,u,v) \,\mbox{ in} \, \Om, \\
		u=v &= 0  \,\mbox{in }\, \Rn \setminus \Om, \\
	\end{aligned}
	\right.
 \end{equation*}
where $f_i$'s are assumed to satisfy the condition $({A1})$, $(A6)$ and $(A8)$. Moreover, we suppose that $h(x)\leq \bar \la$ in $\Om$ and $h(x)<\bar \la$ on a subset of $\Om$ having positive measure, where $\bar \la = \min\limits_{i=1,2}\{\la_{1,p_i}\}$. Now we head to prove Theorem \ref{Thm.1}, so we need the following lemma.
\begin{lemma}\label{lem-1}
	Let $h \in L^{\infty}(\Om)$ with the properties that $ h(x) \leq \bar{\la}$ and
	$h(x)<\bar{\la}$ on a subset of $\Om$ with positive measure. Then, there exists $M>0$ such that
	\begin{equation}\label{Eq:bound}
	\begin{aligned}
	\|u\|_{0,s_1,p_1}^{p_1}+\|v\|_{0,s_2,p_2}^{p_2}-\Iom h(x) (|u|^{p_1}+|v|^{p_2})\,dx
	\geq 2M \bigg( \|u\|_{0,s_1,p_1}^{p_1}+\|v\|_{0,s_2,p_2}^{p_2}    \bigg).
	\end{aligned}
	\end{equation}
\end{lemma}
\begin{proof}
	Suppose (\ref{Eq:bound}) is not true. Then, for each $n \in \N$, there exists
	$(u_n,v_n) \in X$ such that
	\begin{equation*}
		\|u_n\|_{0,s_1,p_1}^{p_1}+	\|v_n\|_{0,s_2,p_2}^{p_2}-\Iom h(x) \big( |u_n|^{p_1}+|v_n|^{p_2}  \big) \, dx<\frac{1}{n} \,\, \text{for all}\, n\in \N,
	\end{equation*}	
	and
	\begin{equation}\label{eq:-1}
	\|u_n\|_{0,s_1,p_1}^{p_1}+\|v_n\|_{0,s_2,p_2}^{p_2}=1.
	\end{equation}
	Therefore, we have,
	\begin{equation}\label{eq:-2}
	\lim_{n \to \infty} \bigg( \|u_n\|_{0,s_1,p_1}^{p_1}+\|v_n\|_{0,s_2,p_2}^{p_2}-\Iom h(x) \big( |u_n|^{p_1}+|v_n|^{p_2}  \big) \,dx \bigg)=0.
	\end{equation}	
	Since $\|u_n\|_{0,s_1,p_1}^{p_1}+\|v_n\|_{0,s_2,p_2}^{p_2}=1$, so $\{u_n\}$ is bounded in $X_{0,s_1,p_1}$ and $\{v_n\}$ is bounded in $X_{0,s_2,p_2}$.
	Therefore, up to a subsequence, we may assume that there exists $u \in X_{0,s_1,p_1}$ and $v \in X_{0,s_2,p_2}$ such that
	\begin{gather*}
		u_n \deb u \,\text{ weakly in}\, X_{0,s_1,p_1} \,\text{and}\, u_n \to u \,\text{ strongly in }\, L^{p_1}(\Om),
	\end{gather*}
	and
	\begin{gather*}
		v_n \deb v \,\text{ weakly in}\, X_{0,s_2,p_2} \,\text{and}\, v_n \to v \,\text{ strongly in }\, L^{p_2}(\Om).
	\end{gather*}
	As $h \in L^{\infty}(\Om)$, we obtain	
	\begin{gather}\label{eq:-3}
		\lim_{n \to \infty} \Iom h(x) \big[ |u_n|^{p_1}+|v_n|^{p_2} \big] \,dx=\Iom h(x) \big[ |u|^{p_1}+|v|^{p_2} \big] \,dx.
	\end{gather}
	Using (\ref{eq:-1}) and (\ref{eq:-3}), from (\ref{eq:-2}), we have,
	\begin{gather}\label{eq:-4}
		1=\lim_{n \to \infty} \big( \|u_n\|_{0,s_1,p_1}^{p_1} +\|v_n\|_{0,s_2,p_2}^{p_2}   \big)=\Iom h(x) \big[ |u|^{p_1}+|v|^{p_2} \big] \,dx.
	\end{gather}
	Since $u_n \deb u$ in $X_{0,s_1,p_1}$ and $v_n \deb v$ in $X_{0,s_2,p_2}$, by the property of weak lower semicontinuity of norm, we get,
	\begin{gather}\label{eq:-5}
		\|u\|_{0,s_1,p_1}^{p_1} +\|v\|_{0,s_2,p_2}^{p_2}
		\leq \liminf_{n \to \infty} \big( \|u_n\|_{0,s_1,p_1}^{p_1} +\|v_n\|_{0,s_2,p_2}^{p_2}   \big).	
	\end{gather}	
	Combining (\ref{eq:-4}), (\ref{eq:-5}) and using the definition of $\la_{1,p_1},\;\la_{1,p_2}$ we have,
	\Bea
	\bar{\la} \big( |u|_{p_1}^{p_1}+|v|_{p_2}^{p_2}  \big) &\leq& \la_{1,p_1} |u|_{p_1}^{p_1}+\la_{1,p_2} |v|_{p_2}^{p_2}\\
	&\leq& \|u\|_{0,s_1,p_1}^{p_1}+\|v\|_{0,s_2,p_2}^{p_2}\\
	&\leq& \Iom h(x) \big( |u|_{p_1}^{p_1} +|v|_{p_2}^{p_2}  \big)\,dx\\
	&\leq& \bar{\la} \big( |u|_{p_1}^{p_1}+|v|_{p_2}^{p_2}  \big).
	\Eea	
	Hence, we have,
	$$
	\Iom \big( \bar{\la}-h(x)   \big) ( |u|_{p_1}^{p_1}+|v|_{p_2}^{p_2}) \, dx=0.
	$$	
	As $h(x)<\bar{\la} $ on a set of positive measure, this yields us,
	$|u|=0=|v|$ a.e. in $\Om$.
	This contradicts (\ref{eq:-4}). This concludes 	our result.
\end{proof}	

\noi \textbf{Proof of Theorem \ref{Thm.1}}. Let us consider the functional $I:X \to \R$ corresponding to $(\mathcal{P}_0)$ defined by
\begin{gather*}
	I(u,v)=\frac{1}{p_1}\|u\|_{0,s_1,p_1}^{p_1}+\frac{1}{p_2}\|v\|_{0,s_2,p_2}^{p_2}-\int_\Om F(x,u,v)dx,\,\text{for all}\, (u,v) \in X.
\end{gather*}
We will show that $I$ has a critical point in $X$ using usual variational technique
which in turn will be a solution of $(\mathcal{P}_0)$.

\noi \textit{\bf Coercivity of $I$}. We note that
\begin{equation}\label{Eq:-I}
\begin{aligned}
I(u,v)
&\geq \frac{1}{p_1 \vee p_2 } \big( \|u\|_{0,s_1,p_1}^{p_1}+ \|v\|_{0,s_2,p_2}^{p_2}  \big)
-\int_\Om F(x,u,v)dx.
\end{aligned}
\end{equation}
By assumption $(A8)$, there exists a function $l \in L^1(\Om)$ such that
\begin{gather}\label{eq:V}
	F(x,u,v) \leq \big( h(x)+\bar{\la} M  \big) \frac{ (|u|^{p_1}+|v|^{p_2}) }{p_1\vee p_2}
	+l(x) \,\text{for all}\, (u,v) \in X
\end{gather}
and $M$ is given in Lemma \ref{lem-1}. Hence, we have,
\begin{equation}\label{eq:bar-V}
\begin{aligned}
\int_\Om F(x,u,v)\,dx
&\leq \Iom h(x) \big( |u|^{p_1} +|v|^{p_2}  \big) + \frac{\bar{\la} M}{p_1 \vee p_2} \Iom ( |u|^{p_1}+|v|^{p_2}) + \Iom l(x)\,dx.
\end{aligned}
\end{equation}
Using (\ref{eq:bar-V}), from (\ref{eq:V}) we obtain  using Lemma \ref{lem-1} that,
\begin{equation}\label{eq:-6}
\begin{aligned}
I(u,v) &\geq \frac{1}{p_1 \vee p_2} \bigg(  \|u\|_{0,s_1,p_1}^{p_1}  +  \|v\|_{0,s_2,q}^{p_2} -\Iom h(x) (|u|^{p_1} +|v|^{p_2})-\bar{\la}M (|u|^{p_1}_{p_1}+|v|_{p_2}^{p_2}) \bigg)-|l|_1  \\
&\geq \frac{1}{p_1 \vee p_2} \bigg(  2M(\|u\|_{0,s_1,p_1}^{p_1}  +  \|v\|_{0,s_2,p_2}^{p_2}) -\bar{\la}M (|u|^{p_1}_{p_1}+|v|_{p_2}^{p_2} ) \bigg)-|l|_1  \\
&\geq \frac{M}{p_1 \vee p_2}
\bigg(  \|u\|_{0,s_1,p_1}^{p_1}  +  \|v\|_{0,s_2,p_2}^{p_2}  ) \bigg)-|l|_1.  \\
\end{aligned}
\end{equation}
Inequality (\ref{eq:-6}) implies $I$ is coercive.

\noi \textit{\bf Weak lower semicontinuity of $I$}. By assumption $(A16)$ we have that 
\begin{equation*}
    \begin{split}
        F(x,u,v) &\leq F(x,0,0)+ \int_0^1(f_1(x,tu,tv)u + f_2(x,tu,tv)v)~dt\\
        & \leq F(x,0,0)+ \frac{C_3}{q}\left(|u|^{q-1}u+ |v|^{q-1}v+  |u|^{q-1}v+ |v|^{q-1}u + u+ v\right)
    \end{split}
\end{equation*}
This implies that for some constant $C>0$ such that 
$$\int_\Om F(x,u,v)~dx \leq C\left(1 + |u|_q +|v|_q+ |u|_1+|v|_1+ \int_\Om(|u|^{q-1}v+ |v|^{q-1}u)~dx\right) $$
{{
Our next aim is to show that $(u,v) \mapsto \Iom F(x,u,v)\,dx$ is weakly lower semicontinuous. Infact, we will show that if $(u_n,v_n) \rightharpoonup (u,v)$ weakly in $X$ then, 		
\begin{equation*}
\lim_{n \to \infty}\Iom F(x,u_n,v_n)\,dx = \Iom F(x,u,v)\, dx.
\end{equation*}	
As $1<q<p_1^* \wedge p_2^*$, by compact embedding, up to a subsequence, we know that		
\begin{gather}\label{eq}
(u_n,v_n) \to (u,v) \;\;\,\text{strongly in}\,\;\; L^q(\Om) \times L^q(\Om).
\end{gather}
Thus, we note that for some constant $C>0$ which may vary at each step,
\begin{equation*}
\begin{aligned}
&\bigg| \Iom F(x,u_n,v_n)\,dx-\Iom F(x,u,v)\, dx\bigg|\\
	&	\leq \Iom \big| F(x,u_n,v_n)- F(x,u,v)\big|\, dx\\
	&\leq \Iom \bigg| \int_0^1 [ F_u(x,tu_n+(1-t)u, tv_n+(1-t)v)(u_n-u) \\
	\quad &+ F_v(x,tu_n+(1-t)u, tv_n+(1-t)v)(v_n-v)] \,dt \bigg|\, dx \\
	&\leq \Iom \bigg| \int_0^1 [ f_1(x,tu_n+(1-t)u, tv_n+(1-t)v)(u_n-u) \\
	&\quad + f_2(x,tu_n+(1-t)u, tv_n+(1-t)v)(v_n-v)] \,dt \bigg|\, dx \\
	&\leq C \Iom \int_0^1 \big( |tu_n+(1-t)u|^{q-1}+|tv_n+(1-t)v|^{q-1}+1  \big) \big( |u_n-u|+|v_n-v|  \big)\,dtdx\\
	&\leq C \Iom \left(|u_n|^{q-1}+|v_n|^{q-1}+|u|^{q-1}+|v|^{q-1}+1\right)(|u_n-u|+|v_n-v|)\, dx\\
	&\leq C \left(\Iom (|u_n|^q+|v_n|^q+|u|^q+|v|^q+1)\, dx  \right)^{\frac{q-1}{q}} \left(\int_\Om  (|u_n-u|^q+|v_n-v|^q)\,dx \right)^{\frac{1}{q}}\\
	&\to 0\,\text{as}\, n \to \infty,
\end{aligned}
\end{equation*}	
where we have used H\"older inequality in the last step, used the fact that $u_n,v_n$ are bounded sequences in $L^q(\Om)$ and $u_n \to u$ and $v_n \to v$ in $L^q(\Om)$. Finally, as norm is weakly lower semicontinuous, so we get that $I$ is weakly lower semicontinuous. From $X$ being a reflexive Banach space, it follows that $I$ has a minimum at some point $(u_0,v_0) \in X$, say. We note that by assumption $(A1)$, $F$ is differentiable. Therefore, $I$ must continuously differentiable on $X$. Since $I$ has a minimum at $(u_0,v_0)$, so we have, $I'(u_0,v_0)=0$.
Hence, it is noted that
\begin{equation*}
	\begin{aligned}
		&\int_{\R^{2N}} \frac{|u_0(x)-u_0(y)|^{p_1-2}(u_0(x)-u_0(y))(w(x)-w(y))}{|x-y|^{N+s_1p_1}}\,dxdy\\
		&+\int_{\R^{2N}} \frac{|v_0(x)-v_0(y)|^{p_2-2}(v_0(x)-v_0(y))(z(x)-z(y))}{|x-y|^{N+s_2p_2}}\,dxdy\\
		&=\Iom u_0f_1(x,u_0,v_0)\,dx+\Iom v_0f_2(x,u_0,v_0)\,dx
	\end{aligned}
\end{equation*}
which suggests that $(u_0,v_0)$ is a solution of $(\mathcal{P})$. This finishes our proof. 

}}

\section{Proof of Theorem \ref{Thm.2}}
{{ This section is devoted to the  proof of Theorem \ref{Thm.2} which goes hand in hand with the proof of Theorem \ref{Thm.1}. For the sake of completeness, we provide the proof.
}}
Let us consider the energy functional $\bar{I}: X \to \R$ defined by
\begin{gather*}
	\bar{I}(u,v)=\bar{J}(u,v)-\int\Om F(x,u,v)dx\,\text{for all}\, (u,v) \in X,
\end{gather*}
where
\begin{equation*}
	\begin{aligned}
		\bar{J}(u,v)&=\int_{\R^{2N}} \phi_1 \big( u(x)-u(y)  \big) (u(x)-u(y)) K_1(x,y)\,dxdy\\
		&+\int_{\R^{2N}} \phi_2 \big( v(x)-v(y)  \big) (v(x)-v(y)) K_2(x,y)\,dxdy,
	\end{aligned}
\end{equation*}
for all $(u,v) \in X$. By the condition {\eqref{*1}} and {\eqref{*2}} defining $\phi_i$ and $K_i, i=1,2$, we note that $\bar{J}$ is convex and lower semicontinuous. By assumption $(A2)$, $F$ is weakly continuous. Hence, $\bar{I}$ is weakly lower semicontinuous. By assumption $(A1)$, there exists a function $l \in L^1(\Om)$ such that
\begin{gather*}
	F(x,u,v) \leq \big( h(x)+\bar{\la}M\big) \frac{(|u|^{p_1}+|v|^{p_2})}{p_1 \vee p_2}+l(x)\, \text{for all}\, (x,u,v) \in \Om \times \R \times \R,
\end{gather*}
and $M$ is given in Lemma \ref{lem-1}.

Hence, we have,
\begin{equation}\label{Eq:bar-V}
\int_\Om F(x,u,v)dx \leq \Iom h(x) \big(\frac{(|u|^{p_1}+|v|^{p_2})}{p_1 \vee p_2}\big)+\frac{\bar{\la}M}{p_1 \vee p_2}\Iom (|u|^{p_1}+|v|^{p_2})+\Iom l(x)\,dx.
\end{equation}
Using \eqref{Eq:bar-V}, \eqref{*1} and {\eqref{*2}}, we obtain, as above,
\begin{equation}\label{Eq:-}
\begin{aligned}
I(u,v)
&\geq \frac{M}{p_1 \vee p_2}
\bigg(  \|u\|_{0,s_1,p_2}^{p_2}  +  \|v\|_{0,s_2,p_2}^{p_2}   \bigg)-|l|_1,  \\
\end{aligned}
\end{equation}
concluding $I$ is coercive. As $F$ and $\phi_i$'s are continuously differentiable, by assumption $(A1)$, $\bar{I'}$ is continuous. Hence, there exists $(u_0,v_0) \in X$ such that
\begin{gather*}
	\bar{I}(u_0,v_0)=\min_{(u,v) \in X} \bar{I}(u,v)\quad\text{ and }\quad\, \bar{I}'(u_0,v_0)=0.
\end{gather*}
 which yields,
\begin{equation}
\begin{aligned}
& \int_{\R^{2N}} \phi_1 \big( u_0(x)-u_0(y)  \big) (w(x)-w(y)) K_1(x,y)\,dxdy\\
&\qquad+\int_{\R^{2N}} \phi_2 \big( v_0(x)-v_0(y)  \big) (z(x)-z(y)) K_2(x,y)\,dxdy\\
&=\Iom f_1(x,u_0,v_0)w(x)\,dx+\Iom f_2(x,u_0,v_0)z(x)\,dx,
\end{aligned}
\end{equation}
for all $w,z \in X$. Therefore, $(u_0,v_0)$ is a solution of $(\mathcal{Q})$.

\section{Acknowledgement}
The first author's research is supported by the Czech Science Foundation, project GJ19--14413Y.

%
%

\end{document}